\documentclass[a4paper, 11pt]{article} 

\usepackage{macro,cite}
\usepackage{epsfig,psfrag}

\begin{document}

\title{{\LARGE
Stabilization under shared communication \\
with message losses and its limitations}}

\author{Hideaki Ishii\\[1mm]
Department of Computational Intelligence
                       and Systems Science\\
Tokyo Institute of Technology\\
4259 Nagatsuta-cho, Midori-ku, Yokohama 226-8502, Japan\\
Email: ishii@dis.titech.ac.jp%
}

\addtolength{\abovedisplayskip}{-.4mm}
\addtolength{\belowdisplayskip}{-.4mm}
\addtolength{\jot}{-.2mm}
\addtolength{\textfloatsep}{-2mm}
\addtolength{\floatsep}{2mm}
\addtolength{\parsep}{12mm}
\setlength{\baselineskip}{3.4ex}

\maketitle

\begin{abstract}
We consider a synthesis problem for a remotely controlled linear system 
where the communication is constrained 
because of the shared and unreliable nature of the channel. 
Modeling the constraints by a periodic transmission scheme and 
random message losses, we present an \Hinfty\ design framework
and study the limitations in the communication 
required for stabilization.
\end{abstract}

%----------------------------------------------------------------------
% PINs:   Ishii   [29153]
% Session code: 

\section{Introduction}
\label{sec:intro}

We consider a remote control system in which the plant
has multiple sensors and actuators connected to a controller over
network channels.  In particular, we follow the approach 
of \cite{Ishii:scl08}
and model two constraints due to the shared and unreliable nature of the channels.
One is a periodic transmission scheme under which
the sensor/actuator nodes take turns to transmit messages in a periodic
manner.
The other constraint is that each transmission is subject to 
random loss or delay due to congestion or error in the communication.
Here, it is assumed that if a message is delayed, then it is considered 
lost.  The losses are modeled as Bernoulli processes where the loss
probabilities are a priori known. 
Further, the controller uses the information regarding the losses 
of the messages that it receives as well as those that it sends;
the latter is realized by the use of acknowledgement messages with 
a one-step delay. 

Under this setup, in \cite{Ishii:scl08}, a synthesis method for 
stochastic stabilization of linear time-invariant plants and for optimization 
under an \Hinfty-type norm criterion has been proposed.
It is based on a necessary and sufficient condition expressed 
in the form of linear matrix inequalities.  Hence, we can employ efficient algorithms
to investigate the effects of the communication constraints 
on control performance.  Another advantage of this approach is that 
because of the norm criterion used in the paper, the design can be viewed as 
a natural extension of deterministic robust control methods. 
We also note that the synthesis method has been found useful
in developing subband coding techniques 
for networked control \cite{IshHar:aut08}.

The focus in this paper is on the limitations in stabilization 
due to the message losses and in particular on upper bounds for 
the loss probabilities above which stabilization cannot be accomplished. 
It turns out that such bounds can analytically be obtained.  
We emphasize that the bounds are expressed in terms of the 
unstable poles of the plant together with parameters in 
the communication scheme.

There are two characteristics of the approach in this paper. 
One is that, by following \cite{SeiSen:05}, 
we view the remote control system as a special case of 
Markovian jump systems (see, e.g., \cite{JiChi:90,CosFra:93}).
This is a natural approach especially when acknowledgements are available
on the controller side.  
%An important consequence is that, as we see in this paper, 
%stability problems can be handled by simple matrix inequalities and
%that the synthesis is convex for fairly general remote control problems.
The other is that we employ the periodic transmission scheme
which has been considered in \cite{LuXieFu:03,IshFra_book:02,HriMor:99}.
In this paper, we show critical bounds on losses which are
generalizations of the results for the simpler case of SISO plant
with single-rate communication.

Similar bounds on critical loss probabilities have 
appeared in the recent literature.  
In the context of remote control, early studies on such probabilities 
include \cite{HadTou:02,Sinopoli:04}.
In \cite{Elia:05}, a synthesis problem for 
stochastic stabilization has been considered, 
and a necessary and sufficient
bound for the state feedback case is found.  
The result is extended in \cite{EliEis:04} to various 
remote control configurations.
One difference from the approach in this paper
is that the controllers are limited to deterministic time-invariant systems.
In \cite{ImeYukBas:06}, it is shown in an LQ type problem 
that the availability of acknowledgement 
messages has a crucial impact on controller designs 
and also on the loss probabilities.
This issue is further studied in \cite{GupHasMur:07,Sinopoli:05}, where
LQG problems over lossy channels are investigated.  
Also, an estimation scheme with filters on the sensor side 
has been proposed in \cite{XuHes:cdc05}.
For the case of nonlinear systems, the treatment of random losses 
in the channel is studied in \cite{MasAbdDor:acc05}.

This paper is organized as follows:  In Section~\ref{sec:prelim},
we introduce a class of stochastic systems and some definitions. 
In Section~\ref{sec:problem}, we review the remote control problem 
considered in \cite{Ishii:scl08}.  Then, the stabilization 
problem of this paper and the main results 
are presented in Section~\ref{sec:siso}.  
To illustrate the results, we provide a numerical example 
in Section~\ref{sec:example}.  Finally, concluding remarks are
given in Section~\ref{sec:concl}.
We note that this paper is an extended version of \cite{Ishii:cdc06}.

\section{Periodic systems with random switchings}
\label{sec:prelim}

In this section, we introduce a class of systems 
called periodic systems with random switchings
and provide some definitions and a preliminary result \cite{Ishii:scl08}.

Consider the following periodic system $G_0$ with random switchings:
\begin{equation}
 \begin{split}
  x_{k+1} &= A_{k,\theta(k)}x_k + B_{k,\theta(k)}w_k,\\
  z_{k}   &= C_{k,\theta(k)}x_k + D_{k,\theta(k)}w_k,
 \end{split}
\label{eqn:PTVsys}
\end{equation}
where $x_k\in\R^n$ is the state,
$w_k\in\R^m$ is the input, $z_k\in\R^p$ is the output, and
$\theta_k\in\Ical_M$ is the mode of the system
with the index set $\Ical_M := \{0,1,\ldots,M-1\}$.
The mode $\theta_k$ is assumed to be an independent and 
identically distributed (i.i.d.) stochastic process
determined by the probabilities
$\alpha_i = \Prob\{\theta_k=i\}$, $i\in\Ical_M$.
The system matrices are $N$-periodic, that is,
e.g., $A_{k+N,i} = A_{k,i}$ for each $i\in\Ical_M$ and $k\in\Z_+$.

Let $\Fcal_{k}$ be the sigma-field generated by $\theta_{[0,k]}$.
We assume that the input $w_k$ is $\Fcal_{k-1}$-measurable for each $k$.
Moreover, $w$ is assumed to be in $l^2$
in the sense that $E[\sum_{k=0}^{\infty}\abs{w_k}^2]<\infty$,
where the expectation $\E[\,\cdot\,]$
is taken over the statistics of $\theta$. Let the norm of such 
signals be $\norm{w}:=E[\sum_{k=0}^{\infty}\abs{w_k}^2]^{1/2}$.
We denote by $\Wcal$ the space of such signals.

For the system $G_0$ in \eqref{eqn:PTVsys}, we employ the following notion 
of stability.  
The system \eqref{eqn:PTVsys} with $w_k\equiv 0$
is said to be \textit{stochastically stable} if
for any initial condition $x_0$,
\[
    \E
     \bigl[
      \sum_{k=0}^{\infty}
        \abs{x_k}^2\bigl|\ x_0
     \bigr]
      < \infty.
\]

%where the expectation $\E[\;\cdot\;]$ is taken over the statistics of $\theta$.
%where $\Theta(k-1) := \{\theta(0),\theta(1),\ldots,\theta(k-1)\}$
%and $\E_{\Theta(k-1)}[\cdot]$ denotes the expectation 
%with respect to $\Theta(k-1)$ ***.
%

We next introduce the $l^2$-induced norm of the system $G_0$.
Suppose that $G_0$ 
is stochastically stable and the initial state is $x_0 = 0$.
Then, we define the $l^2$-induced norm of the system as follows:
\[
   \norm{G_0}
     := \sup_{w\in\Wcal,
               w\neq 0}
              \frac{\norm{z}}{\norm{w}}.
\]

In \cite{Ishii:scl08}, characterizations of the stability
and the norm of the system $G_0$ have been obtained.
They are stated in terms of linear matrix inequalities (LMIs). 
%In particular, necessary and sufficient conditions have been
%derived 
%For later use, 
We present the stability result in the following.

\begin{lemma}\label{lem:MSS}\rm
The system $G_0$ in \eqref{eqn:PTVsys} is stochastically 
stable if and only if there exists 
an $N$-periodic matrix $P_k\in\R^{n\times n}$
such that $P_k=P_k^T>0$ and
\begin{equation}
  \sum_{i\in\Ical_M} \alpha_i A_{k,i}^T P_{k+1} A_{k,i} - P_k < 0~~
        \text{for $k\in\Ical_N$},
\label{eqn:MSS}
\end{equation}
where $\Ical_N:=\{0,\ldots,N-1\}$.
\end{lemma}

\section{Remote control system and its stabilization}
\label{sec:problem}

In this section, we first present the remote control system setup
that has been studied in \cite{Ishii:scl08}.
There, an optimal controller design method under an \Hinfty\ norm
criterion is proposed.  In this paper, we consider the analytic
bounds on the loss probabilities to achieve stabilization of this system. 

\begin{figure}[t]
  \begin{center}
%     \footnotesize
     \small
     \unitlength 0.32mm
     \begin{picture}(100,117)(0,0)

\put(100,109){\vector(-1,0){35}}
\put(65,109){\line(1,0){35}}
\put(82.5,113){\makebox(0,0)[b]{$w$}}
\put(100,109){\line(-1,0){35}}
\put(35,87){\framebox(30,30){$G$}}
\put(35,109){\vector(-1,0){35}}
\put(0,109){\line(1,0){35}}
\put(17.5,113){\makebox(0,0)[b]{$z$}}
\put(35,109){\line(-1,0){35}}
\put(35,95){\line(-1,0){25}}
\put(10,95){\line(1,0){15}}
\put(17.5,99){\makebox(0,0)[b]{$y$}}
\put(25,95){\line(1,0){10}}
\put(35,95){\line(-1,0){25}}
\put(10,95){\vector(0,-1){15}}
\put(0,60){\framebox(20,20){$S_{1}$}}
\put(10,60){\vector(0,-1){15}}
\put(0,25){\framebox(20,20){$\theta_{1}$}}
\put(10,25){\line(0,-1){15}}
\put(10,10){\line(1,0){5}}
\put(15,10){\vector(1,0){20}}
\put(25,14){\makebox(0,0)[b]{$\hat{y}$}}
\put(35,0){\framebox(30,20){$K$}}
\put(65,10){\line(1,0){20}}
\put(75,14){\makebox(0,0)[b]{$\hat{u}$}}
\put(85,10){\line(1,0){5}}
\put(90,10){\vector(0,1){15}}
\put(80,25){\framebox(20,20){$S_{2}$}}
\put(90,45){\vector(0,1){15}}
\put(80,60){\framebox(20,20){$\theta_{2}$}}
\put(90,80){\line(0,1){15}}
\put(90,95){\vector(-1,0){25}}
\put(65,95){\line(1,0){10}}
\put(75,95){\line(1,0){15}}
\put(82.5,99){\makebox(0,0)[b]{$u$}}

\end{picture}
     \caption{Remote control over a shared channel}
     \label{fig:system}
  \end{center}
\vspace*{-2mm}
\end{figure}
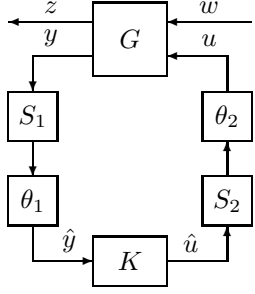

Consider the remote control system depicted in Fig.~\ref{fig:system}.
The generalized plant $G$ is a discrete-time system and
has a state-space equation of the following form:
\begin{equation}
  \begin{split}
     x_{k+1} &= A x_k + B_1 w_k + B_2 u_k,\\
     z_k     &= C_1 x_k + D_{11} w_k + D_{12} u_k,\\
     y_k     &= C_2 x_k + D_{21} w_k, 
  \end{split}
  \label{eqn:G}
\end{equation}
where $x_k\in\R^n$ is the state,
$w_k\in\R^{m_1}$ is the exogenous input, $u_k\in\R^{m_2}$ is 
the control input, $z_k\in\R^{p_1}$ is the controlled output,
and $y_k\in\R^{p_2}$ is the measurement output.
We assume that $(A,B_2)$ is controllable 
and $(A,C_2)$ is observable.

Using a shared communication channel, a remote controller is connected 
to multiple sensors and actuators.  Due to the bandwidth 
limitation in the channel, we assume that at each discrete-time instant,
only one of the sensors or actuators 
can transmit a message over the channel. 
For efficient communication under this constraint, 
the transmission of the messages is periodic
and is a priori fixed.  We now describe this scheme.

Let the period be $N\geq p_2+m_2$.  
%Suppose that there are $p_2$ sensors 
%and that each is capable to transmit messages.  
We index the sensors from 1 to $p_2$ and
fix the order of transmissions within the period $N$ as follows:
Let the index set be $\Ical_{p_2+1}:=\{0,1,\ldots,p_2\}$.
Then, introduce the vector $s_1\in\Ical_{p_2+1}^N$ called
the \textit{switching pattern} for the sensor side.
This specifies that at time $k=lN+r$ with $r=0,1,\ldots,N-1$,
the sensor indexed as $s_{1}(r+1)$ is allowed
to send a message; if $s_{1}(r+1)$ is zero, 
no communication takes place.  
For example, let $N=5$, $p_2=2$, and $s_1 = [1,1,2,0,0]$.
In this case, sensor 1 transmits at $k=0,1,5,6,\ldots$ while
sensor 2 transmits at $k=2,7,\ldots$, and there is no communication
at $k=3,4,8,9,\ldots$.
Similarly, we introduce the switching pattern
$s_2\in\Ical_{m_2+1}^N$ for the $m_2$ actuators; this one determines 
the periodic transmission from the controller to the actuators.
% with the same period $N$.  

We give some notation for the periodic switchings.
First, let $e_{p,i}\in\R^{p}$, $i=1,\ldots,p$, be the unit vectors given by
$e_{p,i}:=[0 \cdots 0\; 1\; 0 \cdots 0]$, where
the $i$th element equals $1$ and the rest are zero.  
%Define two sets $\{s_{1i}\}$ and $\{s_{2i}\}$ of vectors 
%given as follows:
%\begin{align*}
%  s_{1i} 
%    &:= [0 \cdots 0\; 1\; 0 \cdots 0]\in\R^{1\times p_2},~~i=1,\ldots,p_2,\\
%  s_{2j} 
%    &:= [0 \cdots 0\; 1\; 0 \cdots 0]^T\in\R^{m_2\times 1},~~j=1,\ldots,m_2,
%\end{align*}
%where $s_{1i}$ and $s_{2j}$ have the $i$th and $j$th elements equal to 1, 
%respectively, and the rest are zero.  
Now, the \textit{switch boxes} $S_{1}$ and $S_{2}$ in Fig.~\ref{fig:system}
are $N$-periodic matrices defined for $k=lN+r$, $r=0,1,\ldots,N-1$, as
\begin{align*}
  S_{1,k} &:= e_{p_2,i}^T,~~\text{if $s_1(r+1) = i$},\\
  S_{2,k} &:= e_{m_2,j},~~\text{if $s_2(r+1) = j$}.
%  S_{1,k} &:= s_{1i},~~\text{if $s_1(r+1) = i$},\\
%  S_{2,k} &:= s_{2j},~~\text{if $s_2(r+1) = j$}.
\end{align*}

The channel is further constrained by being unreliable due to congestion or delay,
and hence transmitted messages randomly are lost. 
Denote by $\theta_{1,k},\theta_{2,k}\in\{0,1\}$ the stochastic processes
for the message losses, respectively, from the sensors
to the controller and from the controller to the actuators. 
If $\theta_{i,k}=0$, then the message at time $k$ is lost, and otherwise,
it arrives.
They are assumed to be i.i.d.\ Bernoulli processes determined by
\[
  \alpha_1 := \Prob\{\theta_{1,k}=0\}~~\text{and}~~
    \alpha_2 := \Prob\{\theta_{2,k}=0\}
\]
for $k\in\Z_+$.
Letting $\Fcal_{i,k}$ be the sigma-field generated by $\theta_{i,[0,k]}$,
$i=1,2$, we assume 
that the disturbance $w_k$ is $\Fcal_{1,k-1}$- and $\Fcal_{2,k-1}$-measurable.
Further, the disturbance is assumed to be in the
space $\Wcal$ as defined in Section~\ref{sec:prelim}.

The overall plant $\tilde{G}$ including the switches $S_1$ and $S_2$
and the message loss processes $\theta_1$ and $\theta_2$ 
is periodically time varying with period $N$ and with random switchings.  
The state-space equation of $\tilde{G}$ can be expressed as
\begin{equation}
\begin{split}
 x_{k+1}     &= A x_k + B_1 w_k + \theta_{2,k}B_2 S_{2,k} \hat{u}_k,\\
 z_k         &= C_1 x_k + D_{11} w_k + \theta_{2,k} D_{12} S_{2,k}\hat{u}_k,\\
 \hat{y}_k   &= \theta_{1,k}S_{1,k}(C_2 x_k + D_{21} w_k).
\end{split}
\label{eqn:Gtil}
\end{equation}

In Fig.~\ref{fig:system}, $K$ is the controller to be designed.
We allow it to also be $N$-periodic.
Further, we assume that the control $u_k$ is $\Fcal_{1,k}$- and 
$\Fcal_{2,k-1}$-measurable.  
That is, at time $k$, $\theta_{1,k}$ and $\theta_{2,k-1}$ are known 
to the controller.  This is realized by the use of
acknowledgements from the actuators 
regarding the arrival of the control input $u_k$ with a one-step delay.

The controller takes a state-space form as follows:
\begin{equation}
\begin{split}
  \xi_{k}
    &= \hat{A}_{k-1,\theta_1(k-1),\theta_2(k-1)}\xi_{k-1}%\\
%    &\hspace*{2cm}
            + \hat{B}_{k-1,\theta_2(k-1)}\hat{y}_{k-1},\\
  \hat{u}_k 
    &= \hat{C}_{k,\theta_1(k)}\xi_k + \hat{D}_{k}\hat{y}_k,
\end{split}
\label{eqn:K}
\end{equation}
where $\xi_k\in\R^n$ is the state whose dimension is the same as 
that of the plant.  The system matrices are $N$-periodic in $k$: 
For example, $\hat{A}_{k+N,i,j}=\hat{A}_{k,i,j}$ for $k\in\Z_+$ 
and $i,j\in\{0,1\}$. 
Notice that the state equation is expressed for the recursion 
at time $k$.  At this point, $\theta_{2,k-1}$ is available 
at the controller through an acknowledgement. 
Thus, while the $A$- and $B$-matrices make use of this 
information, the $C$- and $D$-matrices can not.
On the other hand, the $B$- and $D$-matrices do not use $\theta_{1,k}$ 
because $\theta_{1,k}=0$ means no input, $\hat{y}_{k} = 0$.

Let the overall closed-loop system in Fig.~\ref{fig:system} be 
$F_l(\tilde{G},K)$.  This system is 
$N$-periodic and has random switchings with 4 modes:
$(\theta_{1,k},\theta_{2,k})=(0,0),(0,1),(1,0),(1,1)$.
It thus falls in the class of systems considered 
in Section~\ref{sec:prelim}.

In \cite{Ishii:scl08}, we have provided an optimal synthesis method 
under an \Hinfty\ criterion.  More specifically, the method
solves the following problem:
For the system in Fig.~\ref{fig:system},
given a scalar $\gamma>0$ and switching patterns $s_1$ and $s_2$,
design a controller $K$ of the form \eqref{eqn:K} such that the 
closed-loop system $F_l(\tilde{G},K)$ is stochastically stable and
satisfies $\norm{F_l(\tilde{G},K)} < \gamma$.
A necessary and sufficient condition for this problem 
has been derived in the form of LMIs.
Thus, using this method, we can numerically check whether 
loss probabilities $\alpha_1$ 
and $\alpha_2$ are small enough to accomplish stabilization.

%--------------------------------------------------------------------

\section{Limitations in loss probabilities for stabilization}
\label{sec:siso}

In this section, we present several upper bounds for the loss probabilities
in the channel which must be met to achieve stabilization.
Hence, the bounds represent the maximum allowable probabilities.
We show that for some specific setups, the bounds become necessary and
sufficient.

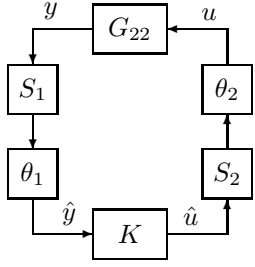
\begin{figure}[t]
  \begin{center}
%     \footnotesize
     \small
     \unitlength 0.32mm
     \begin{picture}(100,105)(0,0)

\put(35,85){\framebox(30,20){$G_{22}$}}
\put(35,95){\line(-1,0){25}}
\put(10,95){\line(1,0){15}}
\put(17.5,99){\makebox(0,0)[b]{$y$}}
\put(25,95){\line(1,0){10}}
\put(35,95){\line(-1,0){25}}
\put(10,95){\vector(0,-1){15}}
\put(0,60){\framebox(20,20){$S_{1}$}}
\put(10,60){\vector(0,-1){15}}
\put(0,25){\framebox(20,20){$\theta_{1}$}}
\put(10,25){\line(0,-1){15}}
\put(10,10){\line(1,0){5}}
\put(15,10){\vector(1,0){20}}
\put(25,14){\makebox(0,0)[b]{$\hat{y}$}}
\put(35,0){\framebox(30,20){$K$}}
\put(65,10){\line(1,0){20}}
\put(75,14){\makebox(0,0)[b]{$\hat{u}$}}
\put(85,10){\line(1,0){5}}
\put(90,10){\vector(0,1){15}}
\put(80,25){\framebox(20,20){$S_{2}$}}
\put(90,45){\vector(0,1){15}}
\put(80,60){\framebox(20,20){$\theta_{2}$}}
\put(90,80){\line(0,1){15}}
\put(90,95){\vector(-1,0){25}}
\put(65,95){\line(1,0){10}}
\put(75,95){\line(1,0){15}}
\put(82.5,99){\makebox(0,0)[b]{$u$}}

\end{picture}
     \caption{Stabilization with periodic transmission}
     \label{fig:systemG22a}
  \end{center}
\vspace*{-2mm}
\end{figure}

For the problem of stabilization, we assume
no disturbance, i.e., $w_k\equiv 0$, throughout this section.  
Hence, we consider the system setup in Fig.~\ref{fig:systemG22a}.
We denote by $G_{22}$ the $(2,2)$-block of the generalized plant $G$. 
To simplify the notation, replace the triple $(A,B_2,C_2)$ with $(A,B,C)$.
Therefore, the realization of $G_{22}$ is given by
\begin{equation}
\begin{split}
   x_{k+1} &= A x_k + B u_k,\\   
   y_k &= C x_k,
\end{split}
\label{eqn:G22}
\end{equation}
where $x_k\in\R^n$, $u_k\in\R^{m_2}$, and $y_k\in\R^{p_2}$.  
We assume that $A$ is an unstable matrix, that is, it has
at least one eigenvalue whose absolute value is larger than $1$.
Denote by $\lambda_i$, $i=0,1,\ldots,n-1$, the eigenvalues of $A$.

%% In this section, we consider a simple setup by assuming the following:
%% (i)~The plant $G_{22}$ is SISO. (ii) The channel is single rate, that is,
%% $S_{i,k}=I$ for $i=1,2$ and $k\in\Z_+$.  Consequently, we deal with
%% the system shown in Fig.~\ref{fig:systemG22b}.
Here, due to the channel on the sensor side, 
the measurement $\hat{y}$ that the controller receives is
\[
   \hat{y}_k = \theta_{1,k} S_{1,k} y_k,
\]
where the behavior of the switch box $S_1$ is determined by
the switching pattern $s_1$.
The control signal $\hat{u}$ transmitted by the controller and
the received signal $u$ at the actuator side are related by
\begin{equation}
   u_k = \theta_{2,k} S_{2,k} \hat{u}_k.
  \label{eqn:uk}
\end{equation}
Similarly, the periodic switch box $S_2$ is specified by 
the switching pattern $s_2$.

In the following, we study three different configurations 
and derive upper bounds for the loss probabilities. 

%% \begin{figure}[t]
%%   \begin{center}
%% %     \footnotesize
%%      \small
%%      \unitlength 0.32mm
%%      \input{./Figures/systemG22b.tex}
%%      \caption{Remote stabilization setup}
%%      \label{fig:systemG22b}
%%   \end{center}
%% \vspace*{-5mm}
%% \end{figure}

\subsection{State feedback under single rate transmission}
\label{sec:state}

The first is the state feedback case with 
a single-rate channel on the actuator side.

In the system in Fig.~\ref{fig:systemG22a},
assume $C=I$ in \eqref{eqn:G22} and $\alpha_1=0$ (that is, 
$\theta_{1,k}\equiv 1$).
The transmissions are assumed to be single rate in that
$S_{i,k}\equiv I$ for $i=1,2$; we thus take $N=1$.
Furthermore, we assume that the plant is single input ($m_2=1$).
Hence, the control input in \eqref{eqn:uk} becomes
\begin{equation}
   u_k = \theta_{2,k} F x_k, 
  \label{eqn:state:u}
\end{equation}
where $F\in\R^{1\times n}$ is the state feedback gain matrix.
Recall that $\theta_{2,k}\in\{0,1\}$ is the loss process
determined by $\alpha_2=\Prob\{\theta_{2,k}=0\}$ for $k\in\Z_+$.

The following proposition provides an upper bound for 
the loss probability $\alpha_2$.  
This result has been given in \cite{Elia:05}.

\begin{proposition}\rm\label{prop:state}
For the system $G_{22}$ in \eqref{eqn:G22}, there exists
a state feedback gain $F$ in \eqref{eqn:state:u} such that
the closed-loop system is stochastically stable 
if and only if the loss probability $\alpha_2$ is sufficiently
small that
\begin{equation}
  \alpha_2 < \frac{1}{\prod_{\abs{\lambda_i}>1} \abs{\lambda_i}^2}.
\label{eqn:state:alpha}
\end{equation}
\end{proposition}

The proof of this proposition is given in the Appendix.
The approach in this paper is based on jump systems results
using Lemma~\ref{lem:MSS}.
We have clarified the critical loss probabilities
in terms of the inequality arising in stochastic stability of such systems 
(see Lemma~\ref{lem:state:1} (ii)).

The result indicates that the unstable poles of the plant
have a direct influence on the allowable loss rates in the channel
for stochastic stability. 
The bound above has been found in \cite{Elia:05}, where
the relation between state feedback stabilization over 
an unreliable channel and an optimal quantizer design
problem has been discussed.  

\subsection{Remote control with periodic transmission}
\label{sec:periodic}

In this subsection, we consider the remote control case with the periodic 
transmission scheme introduced in Section~\ref{sec:problem}.

Consider the system in Fig.~\ref{fig:systemG22a}.
Here, the controller $K$ is limited to the observer-based one as follows:
%similar to \eqref{eqn:remote:obs} in Section~\ref{sec:remote}:
\begin{align}
  \xi_{k} 
    &= A \xi_{k-1} + B u_{k-1} %\notag\\
%    &\hspace*{0.2cm}\mbox{}
         + L_{k-1,\theta_{2}(k-1)}
            \left(
              \theta_{1,k-1} S_{1,k-1} C \xi_{k-1} - \hat{y}_{k-1}
            \right), \label{eqn:periodic:controller}\\
  \hat{u}_k 
    &= F_{k,\theta_{1}(k)}\xi_k, \notag
\end{align}
where $\xi_k\in\R^n$, and $F_{k,\theta_1(k)}$ and
$L_{k,\theta_2(k-1)}$ 
are the feedback and observer gains, respectively. 
These gains are $N$-periodic as, e.g., 
$F_{k+N,i}=F_{k,i}$ for all $k$ and $i$.

We first consider the case when $G_{22}$ is SISO.
Note that in this case, the switch boxes $S_1$ and $S_2$ take values of
either $0$ or $1$ and hence function as discrete-time periodic samplers.
Their periodic behaviors are specified by the switching patterns 
$s_i\in\{0,1\}^N$, $i=1,2$.
We denote by $N_i$ the number of $1$ in the pattern $s_i$ for $i=1,2$.
As a special class of switching patterns, we introduce 
those called \textit{periodic vectors} which take the following form:
\begin{equation}
  s_i = [
          \overbrace{1\,0\,\cdots~0}^{N/N_i}\,\cdots\,
           \overbrace{1\,0\,\cdots\,0}^{N/N_i}
        ]\in\R^{N}.
\label{eqn:s:periodic}
\end{equation}
Some simple examples are $s_i=[1\,0\,\cdots\,0]$ when $N_i=1$, 
and $s_i=[1\,\cdots\,1]$ when $N_i=N$.

We are now ready to state the result on the loss probabilities
for this setup.

%\begin{theorem}\rm\label{thm:periodic}
\begin{proposition}\rm\label{prop:periodic}
  Suppose the switching patterns $s_1$ and 
  $s_2\in\{0,1\}^{N}$ are in the periodic vector form
  in \eqref{eqn:s:periodic}, and have, respectively, 
  $N_1$ and $N_2$ entries of $1$.  
  Then, for the system $G_{22}$ in \eqref{eqn:G22}, 
  there exists a controller $K$ of the form \eqref{eqn:periodic:controller} 
  such that the closed-loop system is stochastically stable if and only if 
  the loss probabilities $\alpha_1$ and $\alpha_2$ satisfy
  \begin{align}
    \alpha_1 
      &< \frac{1}{\prod_{\abs{\lambda_i}>1} \abs{\lambda_i}^{2 N/N_1}},
              \label{eqn:periodic:alpha1}\\
    \alpha_2 
        &< \frac{1}{\prod_{\abs{\lambda_i}>1} \abs{\lambda_i}^{2 N/N_2}}.
              \label{eqn:periodic:alpha2}
  \end{align}
%\end{theorem}
\end{proposition}

\Proofit
Let the estimation error be $e_k := x_k - \xi_k$.
The closed-loop dynamics can be expressed as
\begin{align}
  x_{k+1} &= (A + \theta_{2,k} B S_{2,k} F_{k,\theta_1(k)}) x_k 
                  - \theta_{2,k} B S_{2,k} F_{k,\theta_1(k)} e_k,
   \label{eqn:output:x}\\
  e_{k+1} &= (A + \theta_{1,k} L_{k,\theta_2(k)} S_{1,k} C) e_k.
   \label{eqn:output:e}
\end{align}
The error system \eqref{eqn:output:e} is decoupled from
the state $x_k$ and in particular from the loss process $\theta_{2,k}$.
Using this fact and the structure in the switching pattern $s_1$,
we can show by Lemma~\ref{lem:MSS} that 
to guarantee stochastic stability of this system,
it is necessary and sufficient that there exists 
an observer gain that is independent of $k$ and $\theta_{2,k}$.
The resulting system $e_{k+1} = (A + \theta_{1,k} L S_{1,k} C) e_k$
is periodic with period $\tilde{N}_1:=N/N_1$.

We now look at its dual system  
$z_{k+1} = (A^T + \theta_{1,k} C^T S_{1,k} L^T) z_k$ with
$z_k\in\R^n$ and express it in the so-called 
lifted form as follows: Let $\tilde{z}_k = z_{k\tilde{N}_1}$.  
Then, the lifted system is
\[
  \tilde{z}_{k+1} 
   = \left[
         (A^T)^{\tilde{N}_1} 
            + \theta_{1,k\tilde{N}_1} (A^T)^{\tilde{N}_1-1} C^T L^T
     \right] \tilde{z}_k.
\]
Thus, by Proposition~\ref{prop:state}, 
there exists a stabilizing gain $L^T$ 
if and only if 
\[
  \Prob\{\theta_{1,k\tilde{N}_1}=0\}
   < \frac{1}{%
        \prod_{\abs{\tilde{\lambda}_i}>1}\abs{\tilde{\lambda}_i}^2},
\]
where $\tilde{\lambda}_i$, $i=0,1,\ldots,n-1$, 
are the eigenvalues of $A^{\tilde{N}_1}$.
However, since the loss process $\theta_{1,k}$ is i.i.d.,
the inequality above is equivalent to 
\eqref{eqn:periodic:alpha2}.

Similarly, we can show that the autonomous system 
$x_{k+1} = (A + \theta_{2,k}B S_{2,k} F_{k,\theta_1(k)}) x_k$
of \eqref{eqn:output:x} can be stabilized if and only if
\eqref{eqn:periodic:alpha1} holds.
This implies that a necessary and sufficient 
condition to stabilize the closed-loop system via output feedback
is that the inequalities \eqref{eqn:periodic:alpha1} 
and \eqref{eqn:periodic:alpha2} hold. 
\EndProof

\medskip
This proposition is a generalization of Proposition~\ref{prop:state} 
in two directions.
First, it extends the result for the periodic transmission scheme.
In particular, it shows that the loss probabilities are
constrained by both the unstable dynamics of the plant as well as
the parameters $N$, $N_1$, and $N_2$ in the communication scheme. 
The implication is 
the tradeoff between control performance and transmission rate.
This tradeoff will also be clarified through a numerical 
example in Section~\ref{sec:example}.

Second, the proposition is for the remote control setup with
two communication channels. 
An interesting aspect of the result is that 
the probabilities $\alpha_1$ and $\alpha_2$
for the two channels can be chosen independently and 
further have the same type of maximum for stabilization. 
It can be shown that these characteristics are consequences
of the use of acknowledgement messages; without such
messages, the controller design is no longer convex and 
the analysis becomes much more involved (see, e.g., \cite{Ishii:ifac08}).  

For the single-rate case,
similar problems have been considered in the literature.
In \cite{EliEis:04}, the controller is assumed to be
time invariant, and the approach involves the simultaneous design
of a controller and a decoder on the actuator side; the 
issue of decoder design is considered in the next subsection.
Another work is \cite{Sinopoli:05}, where an LQG problem is studied
for remote control.  

We next consider the case where the plant $G_{22}$ is 
an MIMO system.  The following result is based on 
Proposition~\ref{prop:periodic} and a result in \cite{Ishii:scl08}.

\begin{proposition}\rm\label{thm:periodic_gen}
  Suppose the switching patterns $s_1\in\Ical_{p_2+1}^N$ and 
  $s_2\in\Ical_{m_2+1}^{N}$ have, respectively, 
  $N_1$ and $N_2$ nonzero entries.
  Then, for the system $G_{22}$ in \eqref{eqn:G22},
  there exists a controller $K$ of the form \eqref{eqn:periodic:controller} 
  such that the closed-loop system is stochastically stable only if 
  the loss probabilities $\alpha_1$ and $\alpha_2$ satisfy
  \begin{align}
    \alpha_1 
      &< \frac{1}{\max_{\abs{\lambda_i}>1} \abs{\lambda_i}^{2 N/N_1}},
                 \label{eqn:periodic_gen:alpha1}\\
    \alpha_2 
       &< \frac{1}{\max_{\abs{\lambda_i}>1} \abs{\lambda_i}^{2 N/N_2}}.
                 \label{eqn:periodic_gen:alpha2}
  \end{align}
  The bound on $\alpha_1$ above is also sufficient 
%  if $S_{1,k} \equiv I$ and 
  if $C$ is invertible, and
  the bound on $\alpha_2$ is sufficient 
%  if $S_{2,k} \equiv I$ and 
  if $B$ is invertible.
\end{proposition}

\Proofit
As in the proof of Proposition~\ref{prop:periodic},
the output feedback problem can be separated to 
the state feedback and state estimation problems.  
Hence, we prove only for the state feedback problem.

For this case, we assume $s_1=[1~\cdots~1]$ and $\theta_{1,k}\equiv 1$.
It then follows that the control input is 
\begin{equation}
  u_k=\theta_{2,k} S_{2,k} F_{k,1} x_k.
  \label{eqn:periodic:uk}
\end{equation}
Further, without loss of generality, we assume that $s_2$ has the form
$s_2=[1~\cdots~1~0~\cdots~0]$, where the first $N_2$ entries
are $1$ and the rest are $0$.

We express the system using the lifting technique.
Let $\tilde{x}_k:=x_{kN}$, and let $\tilde{u}$ denote the
lifted signal of $u$:
\[
  \tilde{u}_k 
   := \begin{bmatrix}
       u_{kN}\\ u_{kN+1}\\ \vdots \\ u_{kN+N-1}
     \end{bmatrix},~~k\in\Z_+.
\]
Then, the lifted state equation of the plant can be written as
\[
 \tilde{x}_{k+1} 
   = \tilde{A} \tilde{x}_k + \tilde{B} \tilde{u}_k,
\]
where $\tilde{A} := A^N$ and
$\tilde{B} := \begin{bmatrix}
                   A^{N-1}B & A^{N-2}B & \cdots & B
                \end{bmatrix}$.
Observe that, by the assumption on $s_2$ and by \eqref{eqn:periodic:uk},
we can write $\tilde{u}_k$ as
\begin{align}
  \tilde{u}_k 
    &= \tilde{F}_{\theta_{2}(kN),\cdots,\theta_{2}(kN+N_2-1)}
            \tilde{x}_k,
  \label{eqn:tildeuk}
\end{align}
where
\begin{align*}
   &\tilde{F}_{\theta_{2}(kN),\cdots,
                          \theta_{2}(kN+N_2-1)}
   := \begin{bmatrix}
        \theta_{2,kN} F_{0,1}\\
        \vdots\\
        \theta_{2,kN+N_2-1}F_{N_2-1,1}
           \bigl(
              A+\theta_{2,kN+N_2-2} B F_{N_2-2,1}
            \bigr)\\
        ~~~~~~~~~~~~~~~\cdots
              \bigl(
                 A+\theta_{2,kN}B F_{0,1}
              \bigr)\\
        0\\
        \vdots\\
        0
     \end{bmatrix}.
\end{align*}
Here, we used the facts that $F_{k,i}$ is $N$-periodic in $k$
and $\theta_{1,k}\equiv 1$.
Notice that in the lifted form, the closed-loop system is no longer 
periodic.
In this form, on the other hand, the loss process 
is $(\theta_{2,kN},\cdots,\theta_{2,kN+N_2-1})$, for which
there are $2^{N_2}$ modes for each $k$.

It is clear that the stochastic stability of the original system 
implies that of the lifted system.  Furthermore, it follows that
there exists a gain $\tilde{F}_{1,\ldots,1}$ such that
stability is achieved by the following control:
\begin{align*}
  \tilde{u}_k 
   &= \begin{cases}
       0   %$\tilde{F}_{0,\ldots,0} \tilde{x}_k 
          & \text{if $\theta_{2,kN+i}=0,~\forall i\in\{0,\ldots,N_2-1\}$},\\
       \tilde{F}_{1,\ldots,1} \tilde{x}_k 
          & \text{otherwise}
     \end{cases}
\end{align*}
for $k\in\Z_+$.  
Under this control, the system has only 2 modes, and 
the loss probability is 
\[
   \Prob\left\{
           \text{%
              $\theta_{2,kN+i}=0,~\forall i\in\{0,\ldots,N_2-1\}$}
        \right\} = \alpha_2^{N_2}
\]
for $k\in\Z_+$.
Hence, applying Lemma~\ref{lem:MSS} to this system (with 
$N=1$ in the lemma), we have that
there exist a positive-definite matrix $P\in\R^{n\times n}$ 
and a gain $\tilde{F}_{1,\ldots,1}\in\R^{m_2\times n}$ such that
\begin{align*}
 &\alpha_2^{N_2} \tilde{A}^T P \tilde{A} 
   + (1-\alpha_2^{N_2})
        (\tilde{A}+\tilde{B}\tilde{F}_{1,\ldots,1})^T 
           P (\tilde{A}+\tilde{B}\tilde{F}_{1,\ldots,1})
%  &\hspace*{6.8cm} \mbox{}            
           - P < 0.
\end{align*}
Thus, it follows that 
\[
 \alpha_2^{N_2}
    < \frac{1}{\max_{\abs{\tilde{\lambda}_i}>1} 
                 \abs{\tilde{\lambda}_i}^{2}},
\]
where $\tilde{\lambda}_i$, $i=0,\ldots,n-1$, are 
the eigenvalues of $\tilde{A}$.
Therefore, we arrive at \eqref{eqn:periodic_gen:alpha2}.

To show sufficiency for the case when $B$ is invertible, 
we can take the feedback gain as $F_{k,1}=-B^{-1}A$.
Then, the closed-loop system with the lifted state $\tilde{x}_k=x_{kN}$
is described as
\[
  \tilde{x}_{k+1} 
   = (1-\theta_{2,kN})\cdots(1-\theta_{2,kN+N_2-1})
         A^N \tilde{x}_k
\]
It easily follows from Lemma~\ref{lem:MSS} 
that this system is stochastically stable 
if and only if there exists a matrix $P>0$ satisfying
$\alpha^{N_2} (A^N)^T P A^N - P < 0$.  This is a Lyapunov inequality,
and hence this condition is equivalent to \eqref{eqn:periodic_gen:alpha2}. 
\EndProof

\medskip
We emphasize that in this proposition, 
the switching patterns are not limited to 
those in the periodic vector form as in Proposition~\ref{prop:periodic}.
The bounds however are stated in a very similar form.
On the other hand, in general, the result
is a necessary condition and thus may be conservative.
%The sufficiency part can be further extended 
%to the case when the switching patterns are of periodic-vector type;
%this can be done in a straightforward manner.
It is also remarked that this proposition is a generalization of 
the single-rate version (that is, with $N=1$) 
that has appeared in \cite{ImeYukBas:06,Sinopoli:04,Sinopoli:05,
Katayama:76,Ishii:scl08}; see also \cite{XuHes:cdc05}.

\subsection{Remote control with a decoder}
\label{sec:receiver}

So far, we have assumed that on the actuator side, when
a message is lost, only zero control is applied.
It indeed appears that there might be some room to improve.
In this subsection, we employ a \textit{decoder}
which is a system located at the actuator side to 
compensate the losses as well as the periodic transmission.
It however turns out that for the purpose of stabilization,
the use of such a decoder does not provide advantage. 

\begin{figure}[t]
  \begin{center}
%     \footnotesize
     \small
     \unitlength 0.32mm
     \begin{picture}(130,105)(0,0)

\put(30,85){\framebox(30,20){$G_{22}$}}
\put(30,95){\line(-1,0){20}}
\put(10,95){\line(1,0){20}}
\put(20,99){\makebox(0,0)[b]{$y$}}
\put(30,95){\line(-1,0){20}}
\put(10,95){\vector(0,-1){15}}
\put(0,60){\framebox(20,20){$S_{1}$}}
\put(10,60){\vector(0,-1){15}}
\put(0,25){\framebox(20,20){$\theta_{1}$}}
\put(10,25){\line(0,-1){15}}
\put(10,10){\line(1,0){5}}
\put(15,10){\line(1,0){15}}
\put(30,10){\vector(1,0){20}}
\put(40,14){\makebox(0,0)[b]{$\hat{y}$}}
\put(50,0){\framebox(30,20){$K$}}
\put(80,10){\line(1,0){20}}
\put(90,14){\makebox(0,0)[b]{$\hat{u}$}}
\put(100,10){\line(1,0){15}}
\put(115,10){\line(1,0){5}}
\put(120,10){\vector(0,1){15}}
\put(110,25){\framebox(20,20){$S_{2}$}}
\put(120,45){\vector(0,1){15}}
\put(110,60){\framebox(20,20){$\theta_{2}$}}
\put(120,80){\line(0,1){15}}
\put(120,95){\vector(-1,0){20}}
\put(80,85){\framebox(20,20){$D$}}
\put(80,95){\vector(-1,0){20}}
\put(60,95){\line(1,0){20}}
\put(70,99){\makebox(0,0)[b]{$u$}}

\end{picture}
     \caption{Decoder on the actuator side}
     \label{fig:systemG22r}
  \end{center}
\vspace*{-2mm}
\end{figure}

%For this problem, the setup is as follows:
%We consider the state feedback problem and
%assume only the channel on the actuator side (i.e., $\theta_{1,k}\equiv 1$).
%Extensions to the more general case can be done
%by following the approach in the previous subsections.

Consider the system configuration in Fig.~\ref{fig:systemG22r}.
Again, assume that the plant $G_{22}$ is an SISO system.
Further, we take the switching patterns of both channels in 
the periodic vector form as in \eqref{eqn:s:periodic}.
Recall that in Proposition~\ref{prop:periodic}, the bounds
on the loss probabilities are tight for such cases. 
%We also  restrict the class of controllers $K$ to observer-based ones
%of the form \eqref{eqn:periodic:controller}.

At the actuator side, a decoder $D$ is used.
This is a dynamic, $N$-periodic system that depends on the losses
$\theta_{2,k}$ and outputs the control input $u_k$.
Specifically, it has a state-space form as follows:
\begin{equation}
\begin{split}
  \eta_{k+1} &= A_{D,k} \eta_k + B_{D,k} u_k,\\
  \zeta_k &= C_{D,k} \eta_k,\\
   u_k 
     &= \theta_{2,k} S_{2,k} \hat{u}_k 
           + (1 - \theta_{2,k} S_{2,k}) \zeta_k,
\end{split}
\label{eqn:decoder}
\end{equation}
where $\eta_k\in\R^{n_{D}}$ is the state
and $\zeta_k\in\R$ is the control candidate
produced in the decoder.  The system matrices are $N$-periodic.  
We assume that this system is internally stable in the sense
that if $\theta_{2,k}\equiv 0$, then $\eta_k\rightarrow 0$ 
as $k\rightarrow\infty$.
This guarantees stability in the local feedback of the decoder.

Notice that the control candidate $\zeta_k$ is used when 
a message is lost or when there is no transmission.
A simple example of a decoder is the one-step delay case, where
the decoder functions as a zero-order hold: If a message is not
received, then the previous control value is used.

The following result provides a necessary condition 
for the case with a decoder.

\begin{proposition}\rm\label{prop:decoder}
  Suppose that
  for given switching patterns $s_1$ and $s_2$ in periodic vector forms
  in \eqref{eqn:s:periodic},
  the closed-loop system in Fig.~\ref{fig:systemG22r}
  with the controller $K$ in \eqref{eqn:periodic:controller} 
  and the decoder $D$ in 
  \eqref{eqn:decoder} is stochastically stable.  
  Then, the loss probabilities 
  $\alpha_1$ and $\alpha_2$ satisfy the inequalities in 
  \eqref{eqn:periodic:alpha1}--\eqref{eqn:periodic:alpha2}. 
\end{proposition}

\Proofit
As in the proof of Proposition~\ref{prop:periodic},
the proof can be separated to the state feedback part and
the estimation part.  The estimation part is the same
as in the proposition, and thus the upper bound 
\eqref{eqn:periodic:alpha1} on $\alpha_1$ follows.
We hence show the state feedback part assuming $\theta_{1,k}\equiv 1$.

Letting $\bar{x}_k:=[x_k^T~\eta_k^T]^T$,
we can describe the dynamics of the closed-loop system 
under the state feedback $\hat{u}_{k}=F_{k,\theta_1(k)} x_{k}$ 
as follows:
\[
  \bar{x}_{k+1}
   = \left(
        \bar{A}_k + \theta_{2,k}\bar{B}_k S_{2,k}\bar{F}_k
     \right) \bar{x}_k,
\]
where
\begin{align*}
 \bar{A}_k
  &:= \begin{bmatrix}
       A & B C_{D,k}\\
       0 & A_{D,k} + B_{D,k} C_{D,k}
     \end{bmatrix},~~~
 \bar{B}_k
  := \begin{bmatrix}
           B\\ B_{D,k}
     \end{bmatrix},~~~
 \bar{F}_k
   := \begin{bmatrix}
          F_{k,1} & -C_{D,k}
      \end{bmatrix}.
\end{align*}
Let $\tilde{N}_2:=N/N_2$. Then,
noting that $S_{2,k}=0$ if $k\neq l\tilde{N}_2$, $l\in\Z_+$, we can
lift this system with the lifted state variable 
$\tilde{x}_k:=\bar{x}_{k\tilde{N}_2}$ as
\[
  \tilde{x}_{k+1}
    = \left(
        \tilde{A} + \theta_{2,k\tilde{N}_2} \tilde{B} \tilde{F}
      \right)
        \tilde{x}_{k},
\]
where $\tilde{A} := \bar{A}_{\tilde{N}_2}\cdots\bar{A}_0$, 
$\tilde{B} := \bar{A}_{\tilde{N}_2}\cdots\bar{A}_1 \bar{B}_0$, and 
$\tilde{F} := \bar{F}_0$.
Hence, applying Proposition~\ref{prop:state} to the lifted plant
$(\tilde{A},\tilde{B})$, 
we have that the closed-loop system is stochastically stable only if
\begin{equation}
  \alpha_2 
    < \frac{1}{\prod_{\abs{\tilde{\lambda}_i}>1} \abs{\tilde{\lambda}_i}^2},
 \label{eqn:decoder:alpha2}
\end{equation}
where $\tilde{\lambda}_i$, $i=0,1,\ldots,n+n_{R_0}$, are the eigenvalues
of $\tilde{A}$.  However, the matrix $\tilde{A}$ is upper block diagonal
where the diagonal blocks are $A^{\tilde{N}_2}$ and 
$(A_{D,\tilde{N}_2-1} + B_{D,\tilde{N}_2-1} C_{D,\tilde{N}_2-1})\cdots
(A_{D,0} + B_{D,0} C_{D,0})$.
By the assumption on the decoder, the latter matrix is stable.
Thus, \eqref{eqn:decoder:alpha2} is equivalent to 
\eqref{eqn:periodic:alpha2}.
\EndProof

\medskip
We have several remarks regarding this proposition.
The decoder is an $N$-periodic system and 
can be viewed as a generalized hold device.
It interpolates the control input 
when messages are lost or when no transmission is made.
Clearly, this result and Proposition~\ref{prop:periodic}
imply that, from the perspective of stabilization, 
such decoders are not of help. In fact, it is sufficient to use 
zero control when no message is received by the actuator.

It is however still not clear whether the use of a decoder can improve
the performance of the overall system.  
In the numerical example in the next section,
we make comparisons using an \Hinfty\ design method.
It is also noted that, in general, the design of the decoder $D$
together with the controller $K$ is a difficult problem. 

%% This result is in contrast with \cite{EliEis:04}; there, it is shown
%% that when the plant is nonminimum phase and 
%% when the controller is limited to LTI, 
%% the upper bound \eqref{eqn:remote:alpha} 
%% on the loss probabilities holds true
%% only if a decoder is used.

\section{Numerical example}
\label{sec:example}

We present a numerical example to 
illustrate the results of the paper. 
We consider the system setup in Fig.~\ref{fig:system}
and apply the \Hinfty\ design method introduced in Section~\ref{sec:problem}.

As the generalized plant $G$ in Fig.~\ref{fig:system}, 
we employ the second-order system as follows:
\begin{equation}
\begin{split}
  x_{k+1} &= \begin{bmatrix}
               2   & 0 \\
               0.7 & 1.1
              \end{bmatrix} x_k 
                + \begin{bmatrix}
                     1 \\
                     1
                  \end{bmatrix} w_k 
                    + \begin{bmatrix}
                        1  \\
                        2
                      \end{bmatrix} u_k,\\
  z_k     &= \begin{bmatrix}
               0.5 & -1
             \end{bmatrix} x_k + u_k,\\
  y_k     &= \begin{bmatrix}
               1~& -2
             \end{bmatrix} x_k + w_k.
\end{split}
 \label{eqn:example:plant}
\end{equation}
The system is clearly unstable with eigenvalues $\lambda_i=1.1,2$, $i=1,2$.
We note that the subsystem $G_{22}$ is SISO.  
%Moreover, it is a nonminimum phase system as its transfer function is
%\[
%   G_{22}(z) = \frac{-3 z + 5.5}{z^2 - 3.1 z + 2.2}.
%\]

In the first part of this example, we assumed
perfect transmission on the actuator side
and looked at the effect of the switching
pattern $s_1$ with $N=3$.  Three cases were considered:
$s_1=[1~0~0],[1~1~0],[1~1~1]$.  
According to Proposition~\ref{prop:periodic}, the maximum loss 
probabilities for $\alpha_1$ can be derived for the following 
two cases:
For $s_1=[1~0~0]$, the bound is $1/(\lambda_1\cdot\lambda_2)^6=0.00882$,
and for $s_1=[1~1~1]$, it is $1/(\lambda_1\cdot\lambda_2)^2=0.207$.
The pattern $s_1=[1~1~0]$ is not in the periodic vector form 
\eqref{eqn:s:periodic}, and hence 
the proposition is not applicable. 
We however calculated the probability value
similarly to the one given in the proposition:
$1/(\lambda_1\cdot\lambda_2)^3=0.0939$.

\begin{figure}[t]
  \centering
  \fig{9cm}{6cm}{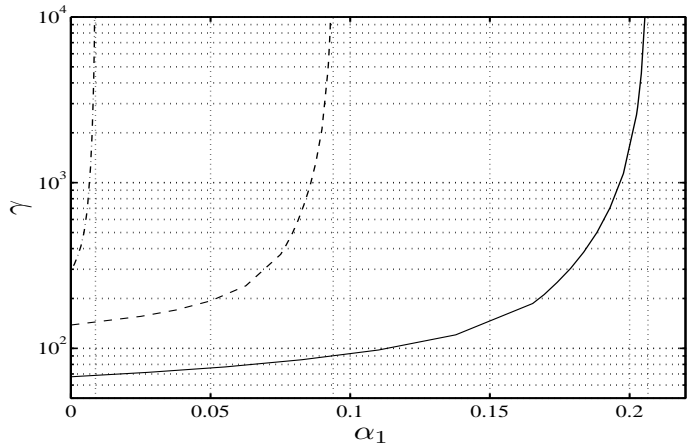}
  \caption{The minimum norm versus $\alpha_1$ for $s_1=[1~0~0]$ (dash-dot),
             $s_1=[1~1~0]$ (dashed), and $s_1=[1~1~1]$ (solid)}
  \label{fig:alpha_gamma_3b}
\vspace*{-1mm}
\end{figure}

A plot showing the minimum \Hinfty\ norm versus $\alpha_1$ 
is given in Fig.~\ref{fig:alpha_gamma_3b}.
It is interesting to note that, for all three cases including $s_1=[1~1~0]$, 
the closed-loop norms explode exactly at the bounds.
%We also notice the fairly large differences in the performance.
This plot exhibits a clear tradeoff between the achievable 
control performance and the transmission rate:
More transmissions at lower loss rate imply better control.

%% \begin{table}[t]
%%  \caption{The maximum loss probabilities}
%%  \centering
%%  \small
%%  \renewcommand{\arraystretch}{1.5}
%%  \label{table}
%%  \begin{tabular}{|c||c|c|c|}
%%     \hline
%%     $s_1$      & $[1~0~0]$ & $[1~1~0]$ & $[1~1~1]$\\
%%     \hline
%%     $\alpha_1$ &  0.00882  &    0.0939 & 0.207\\
%%     \hline
%%  \end{tabular}
%% \vspace*{-3mm}
%% \end{table}

In the second part of simulations,
we assumed a channel only on the actuator side with the
switching pattern $s_2=[1~0]$ and perfect communication on the sensor side.
We designed dynamic controllers of the form \eqref{eqn:K}
for three cases: State feedback, output feedback,
and output feedback with the zero-order hold type decoder
(i.e., it functions as a one-step delay).
In these cases, the upper bound on $\alpha_2$ is 
$1/(\lambda_1\cdot\lambda_2)^4 =0.0427$ 
by the results in Section~\ref{sec:siso}.
We note that the dimension of the controller is different for
the one with the decoder since it was designed for the
generalized plant including the decoder.

%In the state feedback case, 
%we replaced the output equation in \eqref{eqn:example:plant} by
%\[
%   y_k = x_k + \begin{bmatrix}
%                 1\\ 1
%               \end{bmatrix} w_k.
%\]

\begin{figure}[t]
  \centering
  \fig{10cm}{6cm}{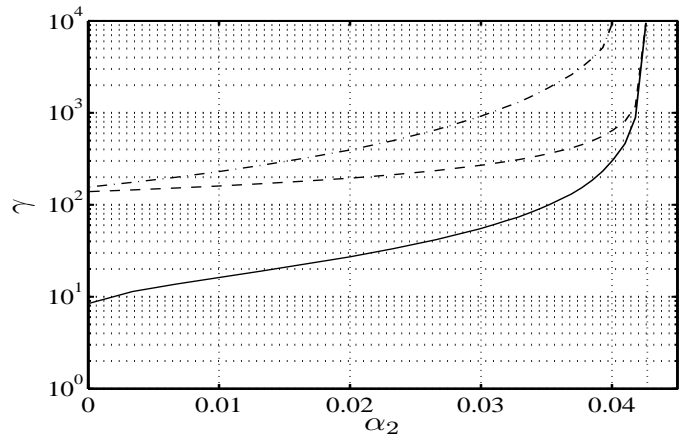}
%  \fig{7.4cm}{5.2cm}{\figdir/alpha_gamma_sep06a.ps}
  \caption{The minimum norm versus $\alpha_2$ for 
           state feedback (solid), output feedback (dashed), and
              output feedback with a decoder (dash-dot)}
  \label{fig:alpha_gamma_decoder}
\vspace*{-1mm}
\end{figure}

For each $\alpha_2$, the minimum \Hinfty\ norm 
for the closed-loop system was calculated.
The results are plotted in Fig.~\ref{fig:alpha_gamma_decoder}.
The norms indeed explode as
$\alpha_2$ approach the upper bound.
It is interesting to note that, for this example, 
the performance of the system with the
decoder is worse especially for large $\alpha_2$.
This may be explained as follows: As $\alpha_2$ 
becomes larger, so does the feedback gain, and hence
the chance to apply a wrong control is higher.

\section{Conclusion}
\label{sec:concl}

In this paper, we have considered the problem of stabilization of 
a linear system
over shared and unreliable channels.  We have shown that there are critical
probability values for the losses above which stability cannot be
achieved.  The implication is the tradeoff between control performance and 
transmission rate for the communication.  
The approach is based on
% results of Markovian jump systems and in particular 
the \Hinfty\ design method proposed in \cite{Ishii:scl08}.
%Future research topics include an \Hinfty\ remote control scheme 
%which does not rely on acknowledgement messages.

\medskip
\noindent
\textit{Acknowledgement}:~
This work was supported in part by
the Ministry of Education, Culture, Sports, Science and Technology, Japan,
under Grant No.\ 17760344.%

%------------------------------------------------------------

{\small
%\bibliographystyle{plain}
%\bibliography{BibDataBase}

}

\appendix
\section{Appendix}

We provide the proof of Proposition~\ref{prop:state}.
This proof consists of two steps.  The next lemma is the key to
relate the state feedback problem to another problem arising in 
\Hinfty\ control.

\begin{lemma}\rm\label{lem:state:1}
The following are equivalent:
\begin{enumerate}
\item[(i)]   The closed-loop system of \eqref{eqn:G22} 
             and \eqref{eqn:state:u} is stochastically stable. 
\item[(ii)]  There exists a positive-definite matrix $P\in\R^{n\times n}$
and a gain $F\in\R^{1\times n}$ such that 
\begin{equation}
    \alpha_2 A^T P A + (1-\alpha_2) (A+BF)^T P (A+BF) - P < 0.
       \label{eqn:state:F}
\end{equation}
\item[(iii)]  There exists a positive-definite matrix $P\in\R^{n\times n}$
             such that 
             \begin{equation}
             \begin{split}
               & A^TPA - P
%               & \hspace*{0.2cm}\mbox{}
                          - A^T P B \left( 
                                        \frac{1}{1-\alpha_2} + B^T P B
                                     \right)^{-1} B^T P A < 0,\\
               & \frac{1}{\alpha_2} - B^T P B > 0.
            \end{split}
            \label{eqn:state:iib}
            \end{equation}
\item[(iv)] There exists a state feedback gain $F$ such that 
\[
   \left\|
     F (zI - A - BF)^{-1} B
   \right\| < \frac{1}{\sqrt{\alpha_2}},
\]
where $\norm{\cdot}$ is the \Hinfty\ norm of a transfer function.
\end{enumerate}
Furthermore, if the condition (iii) holds 
for some $P$, then the inequality \eqref{eqn:state:F}
in (ii) holds with the same $P$ and the gain $F$ given by
\[
   F = -\left( 
           B^T P B 
        \right)^{-1} B^T P A.
\]
\end{lemma}

\Proofit 
The equivalence between (i) and (ii) is a direct consequence of 
Lemma~\ref{lem:MSS}.  
We next show that (ii) holds if and only if (iii) does.  
The inequality \eqref{eqn:state:F} in (ii) is equivalent 
to the following one:
\begin{align}
 &  A^T P A - P       + (1-\alpha_2) \left(
                       F + (B^T P B)^{-1}B^T P A
                     \right)^T  %\notag\\
% & \hspace{2.5cm}\mbox{} \times 
                 B^T P B 
                     \left(
                       F + (B^T P B)^{-1}B^T P A
                     \right)\notag\\
 & \hspace{6cm}\mbox{}  
          - (1-\alpha_2) A^T P B (B^T P B)^{-1} B^T P A < 0. 
\label{eqn:state:PF}
\end{align}
This is shown by expanding \eqref{eqn:state:F} and then completing 
the square for $F$.
Hence, (ii) is equivalent to the existence of $P>0$ such that
\begin{equation}
   A^T P A - P - (1-\alpha_2) A^T P B (B^T P B)^{-1} B^T P A < 0.
\label{eqn:state:pf2}
\end{equation}
Now, suppose that (ii) holds, that is, the inequality \eqref{eqn:state:pf2}
above holds.  This inequality can be expressed as
\begin{align*}
 &  A^T P A - P %\notag\\
% & \hspace{.5cm}\mbox{}
      - A^T P B \left(
                  \frac{\alpha_2}{1-\alpha_2} B^T P B +  B^T P B 
                \right)^{-1} B^T P A < 0.
\end{align*}
We also note that the inequality \eqref{eqn:state:F} holds
for any scaling $\mu P$ of $P$ with positive real $\mu$.
Thus, there exists $\mu$ such that $\mu < (\alpha_2 B^T P B)^{-1}$,
but sufficiently close to $(\alpha_2 B^T P B)^{-1}$ that satisfies
\begin{align*}
  & A^T \mu P A - \mu P %\\
%  &\hspace*{1cm} 
     - A^T \mu P B \left(
                      \frac{1}{1-\alpha_2} + B^T \mu P B
                   \right)^{-1} B^T \mu P A < 0.
\end{align*}
Therefore, for $\mu P$, the inequalities \eqref{eqn:state:iib} 
in (iii) hold.

To show the converse, observe that the second inequality
in \eqref{eqn:state:iib}
implies
\[
  (1-\alpha_2)\left(
                B^T P B
            \right)^{-1} 
   > \left(
       \frac{1}{1-\alpha_2} + B^T P B
     \right)^{-1}.
\]
Hence, for $P$ satisfying \eqref{eqn:state:iib} in (iii),
the inequality \eqref{eqn:state:pf2} also holds.  This implies (ii). 
Moreover, the last statement of the lemma holds true since,
in view of \eqref{eqn:state:PF},
$P$ and the gain $F = -(B^T P B)^{-1} B^T P A$
satisfy the inequality \eqref{eqn:state:F}.  

The equivalence between (iii) and (iv) can be shown 
using the standard \Hinfty\ control theory; see, e.g., \cite{BasBer:95}.
\EndProof

\medskip
The \Hinfty\ control problem in (iii) in the lemma has appeared
in the context of quantized control with 
logarithmic quantizers \cite{EliMit:01,FuXie:05}. 
The following result is from \cite[Lemma~2.4]{FuXie:05}, which
provides an analytic bound for the problem.

\begin{lemma}\rm\label{lem:state:2}
\[
  \inf_{F}
   \left\|
     F (zI - A - BF)^{-1} B
   \right\| = \prod_{\abs{\lambda_i}>1} \abs{\lambda_i},
\]
where the infimum is taken over $F$ such that $A+BF$ is stable.
\end{lemma}

The proof of Proposition~\ref{prop:state} now follows immediately 
from Lemmas~\ref{lem:state:1} and \ref{lem:state:2}.
We note that Lemma~\ref{lem:state:1} provides a method to design 
the feedback gain $F$ for each $\alpha_2$ satisfying 
the bound \eqref{eqn:state:alpha}.

\end{document}